\documentclass{article}%
\usepackage{amsmath}
\usepackage{amsfonts}
\usepackage{amssymb}
\usepackage{graphicx}%
\providecommand{\U}[1]{\protect\rule{.1in}{.1in}}
\newtheorem{theorem}{Theorem}

\newtheorem{corollary}[theorem]{Corollary}

\newtheorem{lemma}[theorem]{Lemma}

\newtheorem{proposition}[theorem]{Proposition}

\newenvironment{proof}[1][Proof]{\noindent\textbf{#1.} }{\ \rule{0.5em}{0.5em}}
\begin{document}

\title{On the finite axiomatizability of some metabelian profinite groups}
\author{Dan Segal}
\maketitle

A group $G$ is \emph{finitely axiomatizable} (\emph{FA}) in a class
$\mathcal{C}$ if $G$ satisfies a first-order sentence $\sigma$ such that every
$\mathcal{C}$ group satisfying $\sigma$ is isomorphic to $G$. Some examples of
this phenomenon where $\mathcal{C}$ consists of profinite groups were
discussed in \cite{NST}; one of the questions raised in that paper was:
\emph{are finitely generated free pro-}$p$ \emph{groups FA}, in either of the
classes profinite groups, pro-$p$ groups?

This is still unknown; a small step in that direction is the following:

\begin{theorem}
\label{main}Each f.g. free metabelian pro-$p$ group on at least two generators
is FA in the class of all profinite groups.
\end{theorem}

The proof depends on

\begin{theorem}
\label{wreath}For each $m,d\geq1$ the profinite wreath product$W_{m,d}%
=\mathbb{Z}_{p}^{(m)}\overline{\wr}\mathbb{Z}_{p}^{(d)}$ is FA in the class of
profinite groups.
\end{theorem}

Here,%
\[
W_{m,d}=\underleftarrow{\lim}_{n\rightarrow\infty}C_{p^{n}}^{(m)}\wr C_{p^{n}%
}^{(d)}.
\]

The analogue in the class of abstract groups of Theorem \ref{main} was
recently established by Kharlampovich, Miasnikov and Sohrabi; they deduce it
from the stronger result that \emph{a free metabelian group is
bi-interpretable with} $\mathbb{Z}$ (\cite{KMS} Theorem 30). The proof of this
is quite elaborate; it seems plausible that the analogue should hold for free
metabelian pro-$p$ groups and $\mathbb{Z}_{p}$, but this seems quite difficult.

Facts from Sections 2 and 5 of \cite{NST} will be used without special
mention. We also often use the fact $\mathcal{H}$:\emph{ every finitely
generated profinite group is Hopfian}, that is, each surjective endomorphism
is an isomorphism.

Logical terms (`formula', `sentence') all refer to the ordinary first-order
language of group theory. \ As discussed in \cite{NST}, `isomorphism' for
profinite groups will mean `continuous isomorphism' (among finitely generated
profinite groups these are actually equivalent, for non-trivial reasons).

\section{A reduction}

A subgroup $H$ of a profinite group $G$ is \emph{definably closed} if there is
a formula $\phi(x)$ such that

(i) for every profinite group $P,$ the subset%
\[
\phi(P):=\{s\in P~\mid~P\models\phi(s)\}
\]
is a closed subgroup, and

(ii) ~$H=\phi(G)$.

\begin{proposition}
\label{basic_case}Let $G$ be a pro-$p$ group. Suppose that $G$ has a definably
closed abelian normal subgroup $A\neq1$ of infinite index such that
$A=\mathrm{C}_{G}(a)$ for each $1\neq a\in A$. Then $G$ satisfies a sentence
$\chi$ such that for any profinite group $H$, if $H\models\chi$ then $H$ is a
pro-$p$ group.
\end{proposition}

For the proof, we combine Lemmas 4.5 and 4.6 of \cite{S} to obtain

\begin{lemma}
\label{pro-pX}Let $\Gamma$ be a profinite group and $A$ a profinite $\Gamma
$-module such that%
\begin{align}
\text{for }a\in A,~x\in\Gamma,~~~ax=a  &  \Longrightarrow(a=0\vee
x=1),\label{hyp}\\
pA+A(\Gamma-1)  &  <A,\label{p-hyp}\\
\bigcap_{1\neq x\in\Gamma}A(x-1)  &  =0. \label{comm}%
\end{align}
Then both $\Gamma$ and $A$ are pro-$p$ groups.
\end{lemma}

Now set $\Gamma=G/A$ in Proposition \ref{basic_case}. Then the conditions
(\ref{hyp}), (\ref{p-hyp}) and (\ref{comm}) are satisfied (see the Remark
following Lemma 4.6 in \cite{S}). So $G$ satisfies a sentence $\alpha$ such
that for any profinite group $H$ satisfying $\alpha$, $H$ has a closed,
definable abelian normal subgroup $B$, and each of (\ref{hyp}), (\ref{p-hyp})
and (\ref{comm}) holds with $B$ for $A$ and $H/B$ for $\Gamma$. It follows by
the preceding lemma that $H$ is a pro-$p$ group.

\section{Ring lemmas}

\begin{lemma}
\label{firstringlemma}Let $R=\mathbb{Z}_{p}[[\zeta]]$ and $M=R^{(m)}$, a free
$R$-module of rank $m\geq1$. Let $K$ be an $R$-submodule of $M$ such that

\begin{itemize}
\item $M/(K+M\zeta)\cong\mathbb{Z}_{p}^{(m)}$

\item $a\zeta\in K\implies a\in K$ for all $a\in M.$
\end{itemize}

Then $K=0$.
\end{lemma}

\begin{proof}
The quotient map $M\rightarrow M/K$ induces an epimorphism%
\[
\mathbb{Z}_{p}^{(m)}\cong M/M\zeta\rightarrow M/(K+M\zeta)\cong\mathbb{Z}%
_{p}^{(m)}.
\]
This must be an isomorphism (by $\mathcal{H}$), so $K\subseteq M\zeta$. Thus
$K=K\cap M\zeta=K\zeta$, which implies that $K=\bigcap_{n\in\mathbb{N}}%
K\zeta^{n}=0$.
\end{proof}

\bigskip

Let $R=\mathbb{Z}_{p}[[\zeta_{1},\ldots,\zeta_{d}]]$. Let us call the tuple
$(\zeta_{1},\ldots,\zeta_{d})$ a \emph{base} for $R$ (together with $p$ it
forms a particular kind of system of parameters for the local ring $R$). The
set of all such bases is denoted $\mathcal{B}(R)$, and for $i<d$ set%
\[
\mathcal{B}_{i}(R)=\left\{  (\zeta_{1},\ldots,\zeta_{i})~\mid(\zeta_{1}%
,\ldots,\zeta_{d})\in\mathcal{B}(R)~\text{\ for some }\zeta_{i+1},\ldots
,\zeta_{d}\right\}  .
\]
Let $\mathcal{X}_{i}\subseteq\mathcal{B}_{i}(R)$ for $1\leq i<d$. I will call
the sequence $(\mathcal{X}_{1},\ldots,\mathcal{X}_{d-1})$ \emph{rich} if
$\mathcal{X}_{1}$ contains infinitely many pairwise non-associate elements,
and for $i>1$ and each fixed $(\zeta_{1},\ldots,\zeta_{i-1})\in\mathcal{X}%
_{i-1}$, there are infnitely many distinct ideals of the form%
\[
\zeta_{1}R+\cdots+\zeta_{i}R
\]
with $(\zeta_{1},\ldots,\zeta_{i})\in\mathcal{X}_{i}$.

\begin{lemma}
\label{rich}Assume that $d\geq2$. For $1\leq i<d$ let $\mathcal{X}_{i}$ be a
subset of $\mathcal{B}_{i}(R)$ such that $(\mathcal{X}_{1},\ldots
,\mathcal{X}_{d-1})$ forms a rich sequence. Then for $1\leq i<d$ we have
\[
D_{i}:=\bigcap_{\overline{\zeta}\in\mathcal{X}_{i}}(\zeta_{1}R+\cdots
+\zeta_{i}R)=0.
\]

\end{lemma}

\begin{proof}
Suppose first that $i=1$. As $\zeta_{1}$ ranges over $\mathcal{X}_{1}$,
$\zeta_{1}R$ ranges over an infinite set of prime ideals of height $1$ in the
Noetherian integral domain $R$, which forces $D_{1}=0$. (This step isn't
really necessary: we could allow $i=1$ in the following argument; but this way
may be less confusing.)

Now let $i>1$ and fix $(\zeta_{1},\ldots,\zeta_{i-1})\in\mathcal{X}_{i-1}$.
Set%
\[
Y=\{\zeta_{i}\in R~\mid~(\zeta_{1},\ldots,\zeta_{i})\in\mathcal{X}_{i}.
\]
Writing $\pi:R\rightarrow\widetilde{R}=R/(\zeta_{1}R+\cdots+\zeta_{i-1}R)$ we
have%
\begin{align*}
D_{i}\pi &  \subseteq\bigcap_{\zeta_{i}\in Y}(\zeta_{1}R+\cdots+\zeta_{i}%
R)\pi\\
&  =\bigcap_{\zeta_{i}\in Y}\widetilde{\zeta}_{i}\widetilde{R}.
\end{align*}
This is the intersection of an infinite set of prime ideals of height $1$ in
the Noetherian integral domain $\widetilde{R}.$ It follows that $D_{i}\pi=0,$
and hence that $D_{i}\subseteq\zeta_{1}R+\cdots+\zeta_{i-1}R$. As $(\zeta
_{1},\ldots,\zeta_{i-1})$ ranges over $\mathcal{X}_{i-1}$ these ideals
intersect in $D_{i-1}$.

The result follows by induction.
\end{proof}

\begin{corollary}
\label{intersectn}Suppose $R=\mathbb{Z}_{p}[[X]]$ where $X$ is the free
abelian pro-$p$ group on $x_{1},\ldots,x_{d}$ and $d\geq2$. Let $\mathcal{C}%
\subseteq X^{(d)}$ denote the set of all bases for $X$. Then%
\[
\bigcap_{\mathbf{y}\in\mathcal{C}}\left(  (y_{1}-1)R+\cdots+(y_{d-1}%
-1)R)\right)  =0.
\]

\end{corollary}

\begin{proof}
The ring $R$ is equal to $\mathbb{Z}_{p}[[\xi_{1},\ldots,\xi_{d}]]$ where
$\xi_{i}=x_{i}-1$ (\cite{DDMS} Thm. 7.20). Set
\[
\mathcal{X}=\left\{  (y_{1}-1,\ldots,y_{d}-1)~\mid~(y_{1},\ldots,y_{d}%
)\in\mathcal{C}\right\}  ,
\]
and let $\pi_{i}:R^{(d)}\rightarrow R^{(i)}$ denote the projection to the
first $i$ factors. Then
\[
\left(  \mathcal{X}\pi_{1},\ldots,\mathcal{X}\pi_{d-1}\right)
\]
is a rich sequence. To see this, note that for $i<d$ and $(y_{1}%
,\ldots,y_{i-1})\in\mathcal{C}\pi_{i-1}$, the group $\widetilde{X}%
=X/\overline{\left\langle y_{1},\ldots,y_{i-1}\right\rangle }$ is free abelian
of rank at least $2,$ and $\mathbb{Z}_{p}[[\widetilde{X}]]$ is naturally
identified with $R/\sum_{j=1}^{i-1}(y_{j}-1)R$. Now $\widetilde{X}$ has
infinitely many $1$-generator direct factors $\overline{\left\langle
\widetilde{y}\right\rangle }$, giving rise to infinitely many distinct
augmentation ideals $(\widetilde{y}-1)\mathbb{Z}_{p}[[\widetilde{X}]]$, the
required condition for a rich sequence. The corollary now follows from the
lemma with $i=d-1$.
\end{proof}

\section{Wreath products\label{wpsec}}

Now we prove Theorem \ref{wreath}.

$W:=W_{m,d}=M\rtimes X$ where $X=\overline{\left\langle x_{1},\ldots
,x_{d}\right\rangle }$ say is a free abelian pro-$p$ group and the $X$-module
$M$ is isomorphic to $R^{(m)}$ where $R=\mathbb{Z}_{p}[[X]]=\mathbb{Z}%
_{p}[[\xi_{1},\ldots,\xi_{d}]]$, writing $\xi_{i}=x_{i}-1$.

Taking $G=W$ and $A=M$ in Proposition \ref{basic_case}, we see that $W$
satisfies a sentence $\chi$ such that every profinite group satisfying $\chi$
is a pro-$p$ group. So it will suffice to show that $W$ is FA in the class of
pro-$p$ groups.

Assume to begin with that $d=1$, and write $x=x_{1}$ etc.

Say $M=a_{1}R\oplus\cdots\oplus a_{m}R$. Then $W$ satisfies a sentence
$\Psi(\mathbf{a},x)$ asserting the following (within the class of pro-$p$ groups):

\begin{itemize}
\item The set $\{a_{1},\ldots,a_{m},x\}$ generates $W$

\item $\overline{\left\langle x\right\rangle }=\mathrm{C}_{W}(x)\cong%
\mathbb{Z}_{p},$

\item $\mathrm{C}_{W}(a_{1})=\ldots=\mathrm{C}_{W}(a_{m}):=M,$ say,

\item $M$ is abelian and normal in $W,$

\item $M/[M,x]\cong\mathbb{Z}_{p}^{(m)},$

\item $M\cap\mathrm{C}_{W}(x)=1$.
\end{itemize}

(In \cite{NST}, \S 5.1 and \S 5.4, it is explained how these are expressed in
first-order language.)

Now suppose that $G$ is a pro-$p$ group and that $G\models\Psi(\mathbf{b},y)$
for some $b_{1},\ldots,b_{m},y\in G$. Write $B$ for the (topological) normal
closure of $\{b_{1},\ldots,b_{m}\}$ and set $Y=\overline{\left\langle
y\right\rangle }$. Then $Y=\mathrm{C}_{G}(y)\cong\mathbb{Z}_{p}$, $B$ is an
abelian normal subgroup contained in $\mathrm{C}_{G}(b_{i})$ for each $i$, and
$G=BY$. It follows that $\mathrm{C}_{G}(b_{i})=B.\mathrm{C}_{Y}(b_{i})=B$ for
each $i$, because $\mathrm{C}_{G}(b_{i})\cap\mathrm{C}_{G}(y)=1$. Thus
$G=B\rtimes Y.$ We consider $B$ as an $R$-module via $x\mapsto y$, and then
$B=b_{1}R+\cdots+b_{m}R$.

Let $K$ be the kernel of the epimorphism $M\rightarrow B$ that sends $a_{i}$
to $b_{i}$ for each $i$. The sentence $\Psi(\mathbf{b},y)$ implies that
$B/B\xi\cong\mathbb{Z}_{p}^{(m)}$ and that $b\xi=0\implies b=0$. It follows
that $K$ satisfies the hypotheses of Lemma \ref{firstringlemma}, and so $K=0$.
Thus $B$ is free of rank $m$ as a module for $\mathbb{Z}_{p}Y,$ and so $G\cong
W$.

Thus $W_{m,1}$ is FA in pro-$p$ groups. Suppose now that $d\geq2$, and
$W=W_{m,d}$. The subgroups $X$ and $M$ are definable by%
\begin{align}
X  &  :=\mathrm{C}_{W}(x_{1})=\ldots=\mathrm{C}_{W}(x_{d})\label{X}\\
M  &  :=\mathrm{C}_{W}(a_{1})=\ldots=\mathrm{C}_{W}(a_{m}). \label{M}%
\end{align}
Let $\Phi(\mathbf{a},\mathbf{x})$ be a first-order formula which asserts (for
the pro-$p$ group $W$) that (\ref{X}) and (\ref{M}) hold and
\begin{align*}
X  &  =\overline{\left\langle x_{1},\ldots,x_{d}\right\rangle }\cong%
\mathbb{Z}_{p}^{(d)},\\
W  &  =\overline{\left\langle a_{1},\ldots,a_{m},x_{1},\ldots,x_{d}%
\right\rangle }\\
\text{ }[M,M]  &  =1,~~M\lhd W
\end{align*}
Suppose that $G$ is a pro-$p$ group and $G\models\Phi(\mathbf{b},\mathbf{y})$
for some $b_{i},$ $y_{j}\in G$. There is an epimorphism $\phi_{\mathbf{b}%
,\mathbf{y}}:W\rightarrow G$ sending $\mathbf{a},\mathbf{x}$ to $\mathbf{b}%
,\mathbf{y}$ respectively. Then $Y=\overline{\left\langle y_{1},\ldots
,y_{d}\right\rangle }\cong\mathbb{Z}_{p}^{d}\cong X$, so $\phi$ induces an
isomorphism from $X$ to $Y$ (in view of $\mathcal{H}$), and so $K_{\mathbf{b}%
,\mathbf{y}}:=\ker\phi_{\mathbf{b},\mathbf{y}}\leq M$.

Suppose that $(t_{1},\ldots,t_{d})$ is a basis for $X$. Denote the
(topological) normal closure of $\{t_{1},\ldots,t_{d-1}\}$ in $W$ by
$N_{\mathbf{t}}$. So%
\[
N_{\mathbf{t}}=[M,t_{1}]\ldots\lbrack M,t_{d-1}]\overline{\left\langle
t_{1},\ldots,t_{d-1}\right\rangle },
\]
and it is easy to see that%
\[
g\in N_{\mathbf{t}}\iff\lbrack W,g]\subseteq\lbrack M,t_{1}]\ldots\lbrack
M,t_{d-1}].
\]
Thus $N_{\mathbf{t}}$ is definable by a formula $\nu(t_{1},\ldots,t_{d}),$ so
by the first case, there is a formula $\Upsilon(t_{1},\ldots,t_{d})$ which
asserts that $W/N_{\mathbf{t}}\cong W_{m,1}$; this statement is true whenever
$(t_{1},\ldots,t_{d})$ is a basis for $X$.

Finally, let $\Theta$ be the sentence asserting, for a pro-$p$ group $G$, that
there exist $b_{1},\ldots,b_{m},y_{1},\ldots y_{d}\in G$ such that (a)
$G\models$ $\Phi(\mathbf{b},\mathbf{y})$ and (b) for each tuple $(s_{1},\ldots
s_{d})$ that generates $Y:=\mathrm{C}_{G}(y_{1})$, $G\models$ $\Upsilon
(\mathbf{s})$.

We have seen that $W$ satifies $\Theta$. Suppose that the pro-$p$ group $G$
satisfies $\Theta$. Then $\phi:=\phi_{\mathbf{b},\mathbf{y}}$ maps $W$ onto
$G$ and $X$ onto $Y$. Let $\mathbf{t}$ be a basis for $X$ and set
$\mathbf{s}=\mathbf{t}\phi$. Then $N_{\mathbf{t}}\phi\leq N_{\mathbf{s}},$ so
we have an induced epimorphism $\phi^{\ast}:F/N_{\mathbf{t}}\rightarrow
G/N_{\mathbf{s}}$. Now $\Upsilon(\mathbf{s})$ asserts that $G/N_{\mathbf{s}%
}\cong W_{m,1}\cong F/N_{\mathbf{t}}$, and it follows by $\mathcal{H}$ that
$K_{\mathbf{b},\mathbf{y}}\leq N_{\mathbf{t}}$. We know that $K_{\mathbf{b}%
,\mathbf{y}}\leq M$, and so
\[
K_{\mathbf{b},\mathbf{y}}\leq M\cap N_{\mathbf{t}}=\sum_{i=1}^{d-1}%
M(t_{i}-1)\text{.}%
\]
Corollary \ref{intersectn} shows that as $\mathbf{t}$ ranges over all bases
for $X$, these modules intersect in zero. It follows that $K_{\mathbf{b}%
,\mathbf{y}}=1,$ and so $G\cong W$.

\section{Free metabelian groups}

$F=F_{d}$ is a free metabelian pro-$p$ group on $d\geq2$ generators
$g_{1},\ldots,g_{d}$. We set $x_{i}=g_{i}F^{\prime}$, $X=F/F^{\prime
}=\overline{\left\langle x_{1},\ldots,x_{d}\right\rangle }$ and $A=F^{\prime}%
$. Then $A$ is a module for the completed group algebra $R=\mathbb{Z}%
_{p}[[X]]$. For $1\leq j\leq d$ set $X_{j}=\overline{\left\langle x_{1}%
,\ldots,x_{j}\right\rangle }$ and $R_{j}=\mathbb{Z}_{p}[[X_{j}]]$.

Note that $R_{j}$ is equal to the power series ring $\mathbb{Z}_{p}[[\xi
_{1},\ldots,\xi_{j}]]$ where $\xi_{i}=x_{i}-1$ for each $i$ (\cite{DDMS} Thm.
7.20). Thus it is a regular local ring of dimension $1+j$.

Write $\Delta_{ji}=\sum_{l=1}^{i}\xi_{l}R_{j};$ the unique maximal ideal of
$R_{j}$ is $\mathfrak{m}_{j}=pR_{j}+\Delta_{jj}$.

Recall that if $G=\overline{\left\langle h_{1},\ldots,h_{d}\right\rangle }$ is
a pro-$p$ group then $G^{\prime}=[h_{1},G]\ldots\lbrack h_{d},G],$ a definably
closed normal subgroup. In particular, $A$ is definably closed in $F$.

We will often use the `Jacobi identity' for metabelian groups,%
\[
\lbrack a,b,c][b,c,a][c,a,b]=1;
\]
this follows at once from the Hall-Witt identity when the derived group is abelian.

Putting $u_{ij}=[g_{i},g_{j}]$ we have

\begin{proposition}
\label{NF}Each element of $A$ is uniquely expressible as%
\begin{equation}
a=\sum_{1\leq i<j\leq d}u_{ij}r_{ij} \label{normalform}%
\end{equation}
with $r_{ij}\in R_{j}$ for each $i$ and $j$.
\end{proposition}

\begin{proof}
The analogue of this result for the abstract free metabelian group, $F_{0}$
say, is established in \cite{MR}, section 6 (cf. also \cite{B}, \cite{BR}).
The \emph{existence} of a representation (\ref{normalform}), also in the
pro-$p$ case, is easily deduced from the Jacobi identity. The proof of
\emph{uniqueness} explained in \cite{MR} uses Fox derivatives; these induce
mappings $d_{j}:F_{0}^{\prime}\rightarrow\mathbb{Z}(F_{0}^{\mathrm{ab}})$
which are $F_{0}$-module homomorphisms and satisfy%
\[
d_{j}(u_{ik})=\left\{
\begin{array}
[c]{ccc}%
0 &  & j\neq i,k\\
x_{k}-1 &  & j=i\\
1-x_{i} &  & j=k
\end{array}
\right.  .
\]
It is easy to verify that each $d_{j}$ extends by continuity to an $R$-module
homomorphism from $A$ to $R$, noting that $A$ is the completion of
$F_{0}^{\prime}$ w.r.t. the $I$-adic topology where $I$ is the ideal of
$\mathbb{Z}(F_{0}^{\mathrm{ab}})$ generated by $p$ and $\xi_{1},\ldots,\xi
_{d},$ while $R$ is the completion of $\mathbb{Z}(F_{0}^{\mathrm{ab}})$ w.r.t.
the $I$-adic topology.

Now suppose that $a$ in (\ref{normalform}) is equal to $0$. We have to show
that each $r_{ij}$ is zero. Arguing by induction on $d,$ we may suppose that
$r_{ij}=0$ for all $i<j<d.$ Then for $1\leq j<d$ we have%
\[
0=d_{j}(a)=(x_{d}-1)r_{jd},
\]
and the result follows.
\end{proof}

\bigskip

The uniqueness of expression in (\ref{normalform}) implies in particular that
$ax_{1}\neq a$ if $a\neq0.$ It follows by symmetry that $ax_{d}\neq a$ if
$a\neq0,$ so the mapping $a\longmapsto a\xi_{d}$ is injective. Noting that for
$i<j<d$ we have $u_{ij}\xi_{d}=u_{id}\xi_{j}-u_{jd}\xi_{i},$ we see that $A$
embeds in a free submodule:

\begin{corollary}
\label{freemod}%
\[
A\cong A\xi_{d}\leq\bigoplus_{i=1}^{d-1}u_{id}R.
\]

\end{corollary}

It follows in turn that $0\neq a\in A$ implies $A=\mathrm{C}_{F}(a)$. So we
may apply Proposition \ref{basic_case} to find a sentence $\chi$, satisfied by
$F$, such that every profinite group satisfying $\chi$ is a pro-$p$ group.
Thus to complete the proof of Theorem \ref{main} it will suffice to show that
$F$ is FA in the class of pro-$p$ groups.

\bigskip Now set
\begin{align*}
C  &  =C_{\mathbf{g}}=\mathrm{C}_{F}(g_{d})\\
H  &  =H_{\mathbf{g}}=AC\\
B  &  =B_{\mathbf{g}}=A\xi_{1}+\cdots+A\xi_{d-1}\\
Z  &  =Z_{\mathbf{g}}=\left\{  a\in A~\mid~a\xi_{d}\in B\right\}  .
\end{align*}
(If $d=2$, this means that $Z=B$.)

It follows from Proposition \ref{NF} that
\[
A/B=Z/B\oplus(D+B)/B
\]
where%
\[
D=\bigoplus_{i=1}^{d-1}u_{id}R,
\]
and that $D\cap B=\bigoplus_{i=1}^{d-1}u_{id}(R\xi_{1}+\cdots+R\xi_{d-1})$.
This implies that
\[
A/Z\cong(D+B)/B\cong D/(D\cap B)\cong S^{(d-1)}%
\]
where $S=\mathbb{Z}_{p}[[\overline{\left\langle x_{d}\right\rangle }]].$

Also $C=\overline{\left\langle g_{d}\right\rangle }$. Thus%
\[
\frac{H}{Z}=\frac{A}{Z}\overline{\left\langle g_{d}\right\rangle }%
\cong\mathbb{Z}_{p}^{(d-1)}\overline{\wr}\mathbb{Z}_{p}=W_{d-1,1},
\]
the wreath product discussed above.

Now all the subgroups mentioned, with the possible exception of $D$, are
definable relative to the parameters $(g_{1},\ldots,g_{d})$. In view of
Theorem \ref{wreath}, there is a formula $\Omega(t_{1},\ldots,t_{d})$ such
that $F\models\Omega(\mathbf{g})$ expresses the fact that $C_{\mathbf{g}%
}=\overline{\left\langle g_{d}\right\rangle }$ and $H_{\mathbf{g}%
}/Z_{\mathbf{g}}\cong W_{d-1,1}$.

Let $\mu$ be a sentence asserting for a pro-$p$ group $G$ that $G$ is
metabelian and that $G/\gamma_{3}(G)\cong F/\gamma_{3}(F)$ (recall that
$F/\gamma_{3}(F)$ is FA in pro-$p$ groups by \cite{NST}, Theorem 5.15). In
particular, if $G\models\mu$ then $G$ is \ generated by $d$ elements.

Suppose now that $G$ is a pro-$p$ group and that $G$ satisfies%
\[
\mu\wedge\left(  \forall t_{1},\ldots,t_{d}\right)  (\beta_{d}(t_{1}%
,\ldots,t_{d})\rightarrow\Omega(t_{1},\ldots,t_{d}))
\]
where $G\models\beta_{d}(t_{1},\ldots,t_{d})$ iff $G=\overline{\left\langle
t_{1},\ldots,t_{d}\right\rangle }$. Let $\theta:F\rightarrow G$ be an
epimorphism and set $K=\ker\theta$.

The induced epimorphism $F/\gamma_{3}(F)\rightarrow G/\gamma_{3}(G)$ is an
isomorphism by $\mathcal{H}$, so $K\leq\gamma_{3}(F)=[A,F]$.

Let $t_{i}=g_{i}\theta$ for each $i$. Then $H_{\mathbf{g}}\theta
=H_{\mathbf{t}}$ and $Z_{\mathbf{g}}\theta\subseteq Z_{\mathbf{t}}$, so
$\theta$ induces an epimorphism $\theta^{\ast}:H/Z\rightarrow$ $H_{\mathbf{t}%
}/Z_{\mathbf{t}}$. Since $t_{1},\ldots,t_{d}$ generate $G$, $G\models
\Omega(\mathbf{t})$, so $H_{\mathbf{t}}/Z_{\mathbf{t}}\cong W_{d-1,1}\cong
H/Z$, whence $\theta^{\ast}$ is an isomorphism (by $\mathcal{H}$); since
$K\leq\lbrack A,F]\leq H$ it follows that $K\leq Z$. Thus%
\[
K\leq Z\cap\lbrack A,F]=Z\cap(B+A\xi_{d})\leq Z\cap D\leq B=B_{\mathbf{g}}.
\]

This holds irrespective of the chosen basis $\mathbf{g}$ for $F$. Now
Corollaries \ref{freemod} and \ref{intersectn} together show that as
$\mathbf{g}$ runs over all such bases, the submodules $B_{\mathbf{g}}$
interesct in zero. Thus $K=1$, and so $G\cong F$.

This completes the proof of Theorem \ref{main}.

\section{The case $d=2$}

There is a much simpler proof when $d=2$. Assume now that $F$ is the free
metabelian pro-$p$ group on generators $g,~h.$ Adapting the notation of the
preceding section, write $A=F^{\prime},$ $x=Ag,$ $y=Ah$, $u=[g,h]$,
$R=\mathbb{Z}_{p}[[F/F^{\prime}]]=\mathbb{Z}_{p}[[\xi,\eta]]$ where $\xi=x-1,$
$\eta=y-1.$ Thus%
\[
A=uR\cong R
\]
by Proposition \ref{NF}, and%
\begin{align*}
\gamma_{3}(F) &  =[A,g][A,h]\\
&  =A\xi+A\eta
\end{align*}
using additive notation for the $F/F^{\prime}$-module $A$.

Set%
\[
H=F/\gamma_{3}(F),
\]
this is the free class-$2$ nilpotent pro-$p$ group (the Heisenberg group over
$\mathbb{Z}_{p}$).

Let $\Psi(s,t)$ be a first-order formula such that for a pro-$p$ group $G$ and
$s,~t\in G$, $G\models\Psi(s,t)$ if and only if

\begin{enumerate}
\item $s$ and $t$ generate $G$

\item $B:=[G,s][G,t]$ is abelian (recall that given \textbf{1.}, $B$ is in
fact the derived group of $G$)

\item $G/[B,s][B,t]\cong H$ (note that given \textbf{2.}, $[B,s][B,t]=\gamma
_{3}(G)$ )

\item For $a,~b\in B$,%
\[
\lbrack a,s][b,t]=1\Longleftrightarrow a=[c,t]\wedge b=[c^{-1},s]\text{ for
some }c\in B.
\]
(Here we use the fact that $H$ is FA in the class of pro-$p$ groups, a special
case of \cite{NST}, Theorem 5.15.)
\end{enumerate}

Now I claim (a) $F$ satisfies $\Psi(g,h)$ and (b) if $G$ is a pro-$p$ group
$G$, $s,~t\in G$, and $G\models\Psi(s,t)$ then $F\cong G$ by a map sending $g$
to $s$ and $h$ to $t$.

This shows that $F$ is FA in the class of pro-$p$ groups; as above we quote
Proposition \ref{basic_case} to infer that $F$ is FA in the class of all
profinite groups.\bigskip

\emph{Proof of }(a). Only condition \textbf{4.} needs comment. Writing $A$
additively, this asserts for $a,~b\in A$ that%
\[
a\xi+b\eta=0\Longleftrightarrow a=c\eta\wedge b=-c\xi\text{ for some }c\in
A\text{.}%
\]
As $A\cong R,$ this follows from the fact that $R\xi\cap R\eta=R\xi\eta$
(while neither of $\xi,$ $\eta$ is a zero-divisor).\bigskip

\emph{Proof of }(b). Now $G$ is a metabelian pro-$p$ group generated by $s$
and $t$, so there exists an epimorphism $\theta:F\rightarrow G$ with
$g\theta=s$ and $h\theta=t$. In view of $\mathcal{H}$, Condition \textbf{3.}
implies that the induced epimorphism $F/\gamma_{3}(F)\rightarrow G/\gamma
_{3}(G)$ is an isomorphism. It follows that $\ker\theta:=K$ is contained in
$\gamma_{3}(F)=[A,g][A,h]$.

Suppose now that $w=[a,g][b,h]\in K$. Then%
\[
1=[a,g][b,h]\theta=[a^{\prime},s][b^{\prime},t]
\]
where $a^{\prime}=a\theta,$ $b^{\prime}=b\theta\in A\theta=B$. According to
\textbf{4.}, there exists $c^{\prime}\in B$ such that $a^{\prime}=[c^{\prime
},t]$ and $b^{\prime}=[c^{\prime-1},s]$. Say $c^{\prime}=c\theta$ for some
$c\in A$ (this exists because $\theta$ maps $A$ onto $B)$. Then%
\[
a=[c,h]w_{1},~~b=[c^{-1},g]w_{2}%
\]
with $w_{1},~w_{2}\in K$. Thus translating into additive notation we have%
\begin{align*}
w &  =[[c,h]w_{1},g][[c^{-1},g]w_{2},h]\\
&  =(c\eta+w_{1})\xi~+~(-c\xi+w_{2})\eta\\
&  =w_{1}\xi+w_{2}\eta.
\end{align*}
It follows that $K\subseteq\lbrack K,F]$. As $F$ is a pro-$p$ group this
forces $K=1$, so $\theta$ is an isomorphism as required.

\end{document}